\documentclass[a4paper,11pt]{article}
\usepackage{times, url}
\usepackage{amsmath,amssymb,amsthm,amsfonts}

\newtheorem{theorem}{Theorem}[section]
\newtheorem{corollary}{Corollary}[theorem]

\newtheorem{remark}[theorem]{Remark}

\begin{document}
\setcounter{page}{1} 
\vspace{10mm}

\begin{center}
{\LARGE \bf  On the Descartes-Frenicle-Sorli and Dris \\
Conjectures Regarding Odd Perfect Numbers}
\vspace{8mm}

{\Large \bf Keneth Adrian Precillas Dagal$^1$ and Jose Arnaldo Bebita Dris$^2$}
\vspace{3mm}

$^1$ Nasser Vocational Training Centre \\ 
Bahrain \\
e-mail: \url{kendee2012@gmail.com}
\vspace{2mm}

$^2$ M.~Sc.~ Graduate, Mathematics Department \\ 
De La Salle University, Manila, Philippines 1004 \\
e-mail: \url{josearnaldobdris@gmail.com}
\vspace{2mm}

\end{center}
\vspace{10mm}

\noindent
{\bf Abstract:} Dris conjectured in his masters thesis that the inequality $q^k < n$ always holds, if $N = {q^k}{n^2}$ is an odd perfect number with special prime $q$. In this note, we initially show that either of the two conditions $n < q^k$ or $\sigma(q)/n < \sigma(n)/q$ holds. This is achieved by first proving that $\sigma(q)/n \neq \sigma(n)/q^k$, where $\sigma(x)$ is the sum of the divisors of $x$.  Using this analysis, we further show that the condition $q < n < q^k$ holds in four out of a total of six cases.  Finally, we prove that $n < q^k$, and that this holds unconditionally.  This finding disproves both the Dris Conjecture and the Descartes-Frenicle-Sorli Conjecture that $k = 1$. 
\vspace{10mm}

\section{Introduction} 

Let $N = {q^k}{n^2}$ be an odd perfect number with special prime $q$ satisfying $q \equiv k \equiv 1 \pmod 4$ and $\gcd(q, n) = 1$.

Therefore, $q \neq n$.  It follows that either $q < n$ or $n < q$.

Dris \cite{Dris} proved that $n < q$ implies the Descartes-Frenicle-Sorli Conjecture that $k = 1$ \cite{Sorli}.  By the contrapositive, $k > 1$ implies that $q < n$.

Acquaah and Konyagin \cite{AcquaahKonyagin} showed that all the prime factors $r$ of $N$ satisfy $r < (3N)^{1/3}$.  In particular, if $k = 1$, then

$$q < (3N)^{1/3} \Longrightarrow q^3 < 3N = 3qn^2 \Longrightarrow q < n\sqrt{3}.$$

Therefore, regardless of the status of the Descartes-Frenicle-Sorli Conjecture, we know that
$$q < n\sqrt{3}$$
must be true.

Brown \cite{Brown} proved that $q < n$ follows from $k = 1$.  (This shows that $q<n$ holds unconditionally.)  Dris \cite{Dris0} and Starni \cite{Starni} also showed $q < n$ using different methods.

Let $\sigma(x)$ be the sum of the divisors of the positive integer $x$.  Let 
$$I(x) = \sigma(x)/x$$
be the abundancy index of $x$.

\section{Main Results}

In this section, we examine this problem:

{\bf Determine the correct ordering for the following quantities:}
$$\frac{\sigma(q)}{n}, \frac{\sigma(n)}{q^k}, \frac{\sigma(q^k)}{n}, \frac{\sigma(n)}{q}$$

Recall the following results:

$$1 < I(qn) = \frac{\sigma(q)}{n}\cdot\frac{\sigma(n)}{q} \leq I({q^k}n) = \frac{\sigma(q^k)}{n}\cdot\frac{\sigma(n)}{q^k} < 2,$$
and
$$n < q^k \Longrightarrow \{k = 1 \Longleftrightarrow n < q\}.$$

In general, since $1 < {q^k}n$ is deficient (being a proper factor of the perfect number $N = {q^k}{n^2}$), then we have
$$\frac{\sigma(q^k)}{n} \neq \frac{\sigma(n)}{q^k}.$$

In a similar vein,
$$\frac{\sigma(q)}{n} \neq \frac{\sigma(n)}{q}.$$

Note that the following implications are true:
$$\frac{\sigma(q)}{n} < \frac{\sigma(n)}{q} \Longrightarrow q < n\sqrt{2},$$
and
$$\frac{\sigma(n)}{q} < \frac{\sigma(q)}{n} \Longrightarrow n < q.$$

We want to show that
$$\frac{\sigma(q)}{n} \neq \frac{\sigma(n)}{q^k}.$$

Suppose to the contrary that
$$\frac{\sigma(q)}{n} = \frac{\sigma(n)}{q^k}.$$

Since $\gcd(q, n) = 1$ and $q$ is prime, we have:
$$\frac{\sigma(q)}{n} = \frac{\sigma(n)}{q^k} \in \mathbb{N}.$$

This means that
$$1 \leq \frac{\sigma(q)}{n} = \frac{\sigma(n)}{q^k}.$$

But since $\sigma(q) = q + 1$ is even while $n$ is odd, we then have:
$$2 \leq \frac{\sigma(q)}{n} = \frac{\sigma(n)}{q^k}.$$

From the inequality
$$2 \leq \frac{\sigma(q)}{n}$$
we get
$$2\cdot\frac{\sigma(n)}{q} \leq \frac{\sigma(q)}{n}\cdot\frac{\sigma(n)}{q} = I(qn) < 2$$
from which we obtain
$$\frac{\sigma(n)}{q} < 1.$$

But then we finally have
$$\frac{\sigma(n)}{q} < 1 < 2 \leq \frac{\sigma(q)}{n} = \frac{\sigma(n)}{q^k} \leq \frac{\sigma(n)}{q},$$
which is a contradiction.

Consequently, we obtain:
$$\frac{\sigma(q)}{n} \neq \frac{\sigma(n)}{q^k}.$$

We now consider two separate cases: \\
Case 1:
$$\frac{\sigma(q)}{n} < \frac{\sigma(n)}{q^k}$$
Since $k \geq 1$, this implies that
$$\frac{\sigma(q)}{n} < \frac{\sigma(n)}{q^k} \leq \frac{\sigma(n)}{q}.$$
Consequently, under Case 1, we have the condition:
$$\frac{\sigma(q)}{n} < \frac{\sigma(n)}{q}.$$
From a previous remark, we know that this implies $q < n\sqrt{2}$. \\
Case 2:
$$\frac{\sigma(n)}{q^k} < \frac{\sigma(q)}{n}$$
Again, since $k \geq 1$, this implies that
$$\frac{\sigma(n)}{q^k} < \frac{\sigma(q)}{n} \leq \frac{\sigma(q^k)}{n}.$$
This implies that, under Case 2, we have the condition:
$$\frac{\sigma(n)}{q^k} < \frac{\sigma(q^k)}{n}.$$
But recall that we have the following inequality \cite{Dris}:
$$\frac{\sigma(q^k)}{q^k} = I(q^k) < \sqrt[3]{2} < I(n) = \frac{\sigma(n)}{n}.$$
Together, the last two inequalities imply that:
$$n < q^k.$$
This implies that the biconditional $k = 1 \Longleftrightarrow n < q$ is true.

We now summarize the results we have obtained so far:
$$\frac{\sigma(q)}{n} \neq \frac{\sigma(n)}{q^k}.$$

The following inequations are trivial:
$$\frac{\sigma(q^k)}{n} \neq \frac{\sigma(n)}{q^k}$$
$$\frac{\sigma(q)}{n} \neq \frac{\sigma(n)}{q}.$$

Also, note that
$$\frac{\sigma(q)}{n} \leq \frac{\sigma(q^k)}{n}$$
and
$$\frac{\sigma(n)}{q^k} \leq \frac{\sigma(n)}{q}.$$

\section{Synopsis}

We now list all the possible orderings for:
$$\left\{\frac{\sigma(q)}{n}, \frac{\sigma(n)}{q^k}, \frac{\sigma(q^k)}{n}, \frac{\sigma(n)}{q}\right\}$$
A:
$$\frac{\sigma(q)}{n} \leq \frac{\sigma(q^k)}{n} < \frac{\sigma(n)}{q^k} \leq \frac{\sigma(n)}{q}$$
B:
$$\frac{\sigma(q)}{n} < \frac{\sigma(n)}{q^k} < \frac{\sigma(q^k)}{n} \leq \frac{\sigma(n)}{q}$$
C:
$$\frac{\sigma(q)}{n} < \frac{\sigma(n)}{q^k} < \frac{\sigma(n)}{q} < \frac{\sigma(q^k)}{n}$$
D:
$$\frac{\sigma(n)}{q^k} < \frac{\sigma(q)}{n} < \frac{\sigma(n)}{q} \leq \frac{\sigma(q^k)}{n}$$
E:
$$\frac{\sigma(n)}{q^k} < \frac{\sigma(q)}{n} < \frac{\sigma(q^k)}{n} < \frac{\sigma(n)}{q}$$
F:
$$\frac{\sigma(n)}{q^k} = \frac{\sigma(n)}{q} < \frac{\sigma(q)}{n} = \frac{\sigma(q^k)}{n}$$

Note that, under cases B, C, D, E and F, we have the inequality
$$\frac{\sigma(n)}{q^k} < \frac{\sigma(q^k)}{n}$$
which implies that $n < q^k$.

Furthermore, note that, under cases A, B, C, D and E, we have the condition
$$\frac{\sigma(q)}{n} < \frac{\sigma(n)}{q}.$$

Lastly, notice that, under cases B, C, D and E, we actually have the inequalities
$$q < n < q^k$$
since $k > 1$ in each of these cases.

\section{A Proof for $n < q^k$}

Let $N = {q^k}{n^2}$ be an odd perfect number with special prime $q$.

It was conjectured in \cite{Dris2} and \cite{Dris} that the inequality $q^k < n$ holds.

Brown \cite{Brown} showed that the Dris Conjecture (that $q^k < n$) holds in many cases.

It is trivial to show that $n^2 - q^k \equiv 0 \pmod 4$.  This means that $n^2 - q^k = 4z$, where it is known that $4z \geq {10}^{375}$.  (See this MSE question \url{https://math.stackexchange.com/q/3556316} and answer \url{https://math.stackexchange.com/a/3557070/28816}, where the case $n < q^k$ is considered.)  Note that if $q^k < n$, then $$n^2 - q^k > n^2 - n = n(n - 1),$$
and that
$${10}^{1500} < N = q^k n^2 < n^3$$
where the lower bound for the magnitude of the odd perfect number $N$ is due to Ochem and Rao (2012) \cite{OchemRao}.  This results in a larger lower bound for $n^2 - q^k$.  Therefore, unconditionally, we have 
$$n^2 - q^k \geq {10}^{375}.$$
We now endeavor to disprove the Dris Conjecture.

Consider the following sample proof arguments:

\begin{theorem}\label{Theorem1}
If $N_1 = q^k n^2$ is an odd perfect number (with special prime $q$) satisfying $n^2 - q^k = 8$, then $n < q^k$.
\end{theorem}

\begin{proof}
Let $N_1 = q^k n^2$ be an odd perfect number (with special prime $q$) satisfying $n^2 - q^k = 8$.

Then $$(n + 3)(n - 3) = n^2 - 9 = q^k - 1.$$

This implies that $(n + 3) \mid (q^k - 1)$, from which it follows that
$$n < n + 3 \leq q^k - 1 < q^k.$$
We therefore conclude that $n < q^k$.
\end{proof}

\begin{theorem}\label{Theorem2}
If $N_2 = q^k n^2$ is an odd perfect number (with special prime $q$) satisfying $n^2 - q^k = 40$, then $n < q^k$.
\end{theorem}

\begin{proof}
Let $N_2 = q^k n^2$ be an odd perfect number (with special prime $q$) satisfying $n^2 - q^k = 40$.

Then $$(n+7)(n-7)=n^2 - 49=q^k - 9,$$
from which it follows that
$$(n+7) \mid (q^k - 9)$$
which implies that
$$n < n+7 \leq q^k - 9 < q^k.$$
\end{proof}

Note that, in the proof of Theorem \ref{Theorem2}, $49$ is not the \emph{nearest square} to $40$ ($36$ is), but rather the nearest square \emph{larger} than $40$.

With this minor adjustment in the logic, we would expect the general proof argument to work.

(Additionally, note that it is known that $n^2 - q^k$ is \emph{not a square}, if $q^k n^2$ is an odd perfect number with special prime $q$.  See this MSE question \url{https://math.stackexchange.com/q/3121498} and the answer contained therein \url{https://math.stackexchange.com/a/3122247/28816}.  Alternatively, you may refer to the preprint \cite{DrisSanDiego}.)

So now consider the equation $n^2 - q^k = 4z$.  Note that $q^k < (2n^2)/3$ by \cite{Dris}.  Following our proof strategy, we have:

Subtracting the smallest square that is larger than $n^2 - q^k$, we obtain

$$n^2 - \bigg(\lceil{\sqrt{n^2 - q^k}}\rceil\bigg)^2 = q^k + \Bigg(4z - \bigg(\lceil{\sqrt{n^2 - q^k}}\rceil\bigg)^2\Bigg).$$

Note that it is always the case that
$$\Bigg((n^2 - q^k) - \bigg(\lceil{\sqrt{n^2 - q^k}}\rceil\bigg)^2\Bigg) < 0,$$
if $N = q^k n^2$ is an odd perfect number with special prime $q$.

More so, note that we always have
$$0 < n - \lceil{\sqrt{n^2 - q^k}}\rceil.$$

Consequently, it follows that
$$\Bigg(n + \lceil{\sqrt{n^2 - q^k}}\rceil\Bigg)\Bigg(n - \lceil{\sqrt{n^2 - q^k}}\rceil\Bigg) = q^k - y$$
for some positive integer $y$.  This would imply that
$$\Bigg(n + \lceil{\sqrt{n^2 - q^k}}\rceil\Bigg) \mid (q^k - y)$$
from which we finally have
$$n < \Bigg(n + \lceil{\sqrt{n^2 - q^k}}\rceil\Bigg) \leq q^k - y < q^k.$$

Since all odd perfect numbers satisfy the equation $n^2 - q^k = 4z$, we now have the following result.

\begin{theorem}\label{MainTheorem}
If $N = q^k n^2$ is an odd perfect number with special prime $q$, then the inequality $n < q^k$ holds.
\end{theorem}

We also get the following corollaries.

\begin{corollary}\label{Corollary1}
If $N = q^k n^2$ is an odd perfect number with special prime $q$, then $k \geq 5$.
\end{corollary}

\begin{proof}
Let $N = q^k n^2$ be an odd perfect number with special prime $q$.

By the earlier result of Brown, Dris and Starni mentioned in the Introduction, we have the inequality $q < n$.

Together with the inequality in Theorem \ref{MainTheorem}, this means that
$$q < n < q^k$$
from which it follows that
$$k > 1.$$
This last inequality implies that $k \geq 5$, since it is known that $k \equiv 1 \pmod 4$.  In particular, we now know that the Descartes-Frenicle-Sorli Conjecture that $k=1$ is false.
\end{proof}

\begin{corollary}\label{Corollary2}
If $N = q^k n^2$ is an odd perfect number with special prime $q$, then $q^k > {10}^{500}$.
\end{corollary}

\begin{proof}
By the result in Theorem \ref{MainTheorem}, we have $n < q^k$.  Squaring both sides of the inequality, we get $n^2 < q^{2k}$.  Multiplying both sides by $q^k$, we obtain $N = q^k n^2 < q^{3k}$.  But $N > {10}^{1500}$ by \cite{OchemRao}. It follows that $q^k > {10}^{500}$.
\end{proof}

\begin{remark}
Note that $q^k > {10}^{500}$ significantly improves on previously known results for components of odd perfect numbers.
\end{remark}

\section{Open Problems}

The best course of action at this point would be to come up with a proof for the Descartes-Frenicle-Sorli Conjecture, in order to finally show that there are indeed no odd perfect numbers.  Or perhaps this research direction is misguided, as the recent article \cite{Nielsen} would show?

We leave this problem for resolution by other researchers.

\section*{Acknowledgements} 

The authors would like to thank the anonymous referees for their valuable feedback and suggestions which helped in improving the quality of this manuscript.
 
\makeatletter
\renewcommand{\@biblabel}[1]{[#1]\hfill}
\makeatother

\end{document}